\documentclass[12pt]{paper}

\usepackage{amssymb,amsfonts,amsthm,amsmath}

\usepackage{latexsym,multicol,graphicx}

\usepackage{graphicx, verbatim}

\usepackage{xcolor}
\usepackage{verbatim}

\usepackage{pict2e}

\usepackage[a4paper, hmargin=2.8cm, vmargin={3.5cm,3.5cm}]{geometry}

\usepackage{fancyhdr}
\usepackage[title,titletoc,toc]{appendix}

\newtheorem{theorem}{Theorem}
\newtheorem{lemma}{Lemma}

\newtheorem{corollary}[theorem]{Corollary}

\renewcommand{\phi}{\varphi}

\renewcommand{\emptyset}{\varnothing}

\newcommand{\e}{\varepsilon}

\newcommand{\la}{\lambda}

\newcommand{\bE}{\mathbb E}

\newcommand{\spt}{{\rm spt}}

\def\done{{1\hskip-2.5pt{\rm l}}}

\renewcommand{\le}{\leqslant}
\renewcommand{\ge}{\geqslant}

\newcommand{\bR}{\mathbb R}
\newcommand{\bC}{\mathbb C}
\newcommand{\bZ}{\mathbb Z}
\newcommand{\bD}{\mathbb D}
\newcommand{\bT}{\mathbb T}
\newcommand{\bN}{\mathbb N}
\newcommand{\bQ}{\mathbb Q}

\newcommand{\bP}{\mathbb P}

\begin{document}

\title{Spectra of stationary processes on $\bZ$}

\author{Alexander Borichev,
Mikhail Sodin\thanks{Supported in part by ERC Advanced Grant~692616 and ISF Grant~382/15}, Benjamin Weiss}

\begin{comment}
\author{Alexander Borichev \\

\and
Mikhail Sodin
\thanks{}
\\
School of Mathematical Sciences, Tel Aviv University\\
Tel Aviv 69978, Israel \\
{\tt sodin@post.tau.ac.il}
\and
Benjamin Weiss
\\
Einstein Institute of Mathematics, Hebrew University of Jerusalem \\
Jerusalem 91904 Israel \\
{\tt weiss@math.huji.ac.il}
}
\end{comment}

\maketitle

\noindent\today

\begin{abstract}
We will discuss a somewhat striking spectral property of finitely valued
stationary processes on $\bZ$ that says that if the spectral measure of the
process has a gap then the process is periodic.  We will give some extensions
of this result and raise several related questions.
\end{abstract}

\section{Introduction}\label{sect:intro}

One of the basic features of a complex valued stationary random process
$\xi\colon \bZ\to \bC$ is its {\em spectral measure} $\rho$ and its closed support,
$\spt(\rho)$, called {\em the spectrum} of the process $\xi$. As is well known, the spectral measure
is obtained in the following way. Assuming for simplicity of notation that
$\xi$ has zero mean, one first forms the covariance function
$r(m)=\bE[\xi(0)\bar\xi(m)]$, and then observing that it is positive definite, defines
the measure $\rho$ on the unit circle $\bT$  by $r(m)=\widehat{\rho}(m)$, i.e.,
$r$ is the Fourier transform of $\rho$. Any positive measure on $\bT$ is a spectral measure of a stationary process obtained, for instance, by a familiar Gaussian construction.
The case of finitely valued stationary processes is strikingly different.
We will show that {\em if $\xi$ is a finitely valued stationary process on $\bZ$ and
the spectrum of $\xi$ is not all of $\bT$, then $\xi$ is periodic and therefore
its spectrum is contained in $\{t\colon t^N=1\}$, where $N$ is a period of $\xi$}.

This result bears a close resemblance to a theorem of Szeg\H{o}.
Szeg\H{o}'s original version~\cite[Section~83]{Dienes}
says that if a sequence of Taylor coefficients $f_n$ attains finitely many values and
the sum of the Taylor series $f(z)=\sum_{n\ge 0} f_n z^n$ can be analytically continued through an arc
of the unit circle, then the sequence $(f_n)$ is eventually periodic and $f$ is a rational
function with poles at the roots of unity. Helson's harmonic analysis version of Szeg\H{o}'s theorem~\cite[Section~6.4]{Helson} says that if a sequence $(f_n)_{n\in\bZ}$ attains finitely many values and its Beurling's spectrum (defined below in Section~\ref{subsect:Beurling})
is not all of $\bT$, then the sequence $(f_n)$ is periodic.

We will give two proofs of our result. They are short, not too far from each other, both
use some ideas of Szeg\H{o} and Helson, and each of these proofs will give us
more than we have stated above.

For our first proof we will show the following general result. If $\xi$ is a  finite valued stationary ergodic process, the Beurling spectrum of almost every realization $\xi(n)$
coincides with the support of the spectral measure of the process. In addition we will give a new proof of Helson's result that gives an effective estimate
for the size of the period in terms of the set of values of the process and the size of the gap in the spectrum.

The second proof exploits a relation to the polynomial
approximation problem in the space $L^2(\rho)$. It allows us
to replace the assumption that $ \spt(\rho) \ne \bT $
by a much weaker condition, which in particular forbids
$ \rho $ to have sufficiently deep exponential zeroes.

It would be very interesting to reveal more restrictions on the spectral measures of finitely valued stationary sequences. In this connection we mention that McMillan~\cite{McMillan}
described covariances of stationary processes that attain the values $\pm 1$, the proofs appeared later in the work of Shepp~\cite{Shepp}. Unfortunately, it seems to be very
difficult to apply this description. On the other hand, in the book by
M.~Queff\'elec~\cite{Que} one can find a Zoo of various spectral measures of finitely valued stationary processes of dynamical origin.

Another intriguing question is to extend our main results to stationary processes on
$\bZ^d$ with $d\ge 2$.

In the last section of this note we will turn to unimodular stationary processes on $\bZ$.
Using a result of Eremenko and Ostrovskii~\cite{EO}, we will show that if $\xi$ is
an ergodic process that takes the values from the unit circle and if the spectrum of $\xi$ is
contained in an arc of length less than $\pi$, then $\xi (n) = t s^n$, $n\in\bZ$, with
a unimodular constant $s$ and a random variable $t$ which has
either a uniform distribution on the circle if $s$ is an irrational rotation, or a uniform distribution on a cyclic group
if $s$ is a rational rotation.

\medskip
\centerline{*\quad*\quad*}

We sadly dedicate this note to the memory of Victor Havin, a wonderful person, teacher,
and mathematician. The interplay between harmonic and complex analysis and probability
theory presented here, probably, would be close to his heart.

We thank Fedor Nazarov for several illuminating discussions of the material presented here, and especially for suggestion
to use monic polynomials with small $ L^2 (\rho) $ norms.

\section{Spectra of stationary random sequences}\label{sect:spectra}

Some of our results deal with wide-sense stationary processes.
Recall that the random process $\xi\colon \bZ\to\bC$ is called {\em wide-sense stationary}
(a.k.a. weak stationary or second order stationary) if $\bE|\xi (n)|^2 <\infty$
for every $n$ and $\bE\, \xi (n)$ and $\bE[\xi(n)\bar\xi(n+m)]$ do not depend on $n$.
Then the spectral measure $\rho$ of $\xi$ is defined in the same way as above so that
the covariance function of $\xi$ is the Fourier transform of $\rho$. As above, we will
call the closed support of $\rho$, $\spt(\rho)$, {\em the spectrum} of $\xi$.

\subsection{Spectrum of a single realization}\label{subsect:spectr_realization}

Now, we define the spectrum of a single realization of a wide-sense stationary
process $\xi$. First, we note that by the Borel-Cantelli lemma, almost surely, we have
\[
|\xi (n)| = o(|n|^\alpha), \quad n\to\infty\,,
\]
with any $\alpha>\tfrac12$. Hence, the sequence $(\xi(n))_{n\in\bZ}$ generates the
random distribution (or the random functional) $F_\xi$ on $C^\infty(\bT)$ acting as
follows:
\[
F_\xi (\phi) = \sum_{n\in\bZ} \xi (n) \widehat{\phi}(-n), \quad \phi\in C^\infty (\bT)\,.
\]
We denote by $\sigma (\xi)\subset\bT$
the support of this distribution, that is, the complement to the largest open set
$O\subset\bT$ on which $F_\xi$ vanishes (which means that $F_\xi (\phi)=0$ for any
smooth function $\phi$ whose closed support is contained in $O$).

We consider the Borel structure on the family of the compact subsets of the unit circle corresponding to the Hausdorff distance.
Then the map $\xi\mapsto \sigma (\xi)$ is measurable.
Indeed, we can take a countable collection of smooth functions
$\{\phi_k\}$ such that the complement of the support of
$\sigma(\xi)$ is determined by the set of $k$
such that $F_{\xi}(\phi_k)=0$. With a natural Borel structure
on the set of sequences that grow polynomially
this is a measurable set for a fixed $k$. Finally, one can characterise the Hausdorff distance
between two compact subsets of the unit circle in terms of the open arcs that
constitute theirs complements.

Note that the random compact set $\sigma(\xi)$ is invariant with respect to the
translations of $\xi$, so if the translations act ergodically on $\xi$, the spectrum
$\sigma (\xi)$ is not random. We also note that $\sigma (\xi)$ is invariant with
respect to the flip $\xi(n)\mapsto \xi(-n)$.

\subsection{Spectrum of the process and spectra of its realizations}
\label{subsect:relation-between-spectra}

\begin{theorem}\label{thm:spectra}
Suppose $\xi\colon\bZ\to\bC$ is a wide-sense stationary process with zero mean.
Let $\rho$ be the spectral measure of $\xi$. Then,

\smallskip\noindent{\rm (A)} almost surely, $\sigma(\xi)\subseteq \spt(\rho)$;

\smallskip\noindent{\rm (B)} if $O\subsetneqq\bT$ is an open set such that,
almost surely, $O\cap\sigma(\xi)=\emptyset$, then $O\cap\spt(\rho)=\emptyset$.
\end{theorem}

\begin{corollary}\label{cor:spectra}
Let $\xi$ be a stationary square integrable ergodic process on $\bZ$. Then, almost surely,
$\sigma (\xi)=\spt(\rho)$.
\end{corollary}

\subsubsection{Proof of Part (A)}
Suppose that $\spt(\rho)\ne\bT$ and take an arbitrary function
$\phi\in C_0^\infty (\bT\setminus\spt(\rho))$. Then
\begin{equation}\label{eq:vanishing}
\int_{\bT}\, \Bigl| \sum_{n\in\bZ} \widehat{\phi}(n) t^n \Bigr|^2\, {\rm d}\rho (t) = 0\,.
\end{equation}
Recall that there exists a linear isomorphism between the closure of the linear span
of $(\xi (n))_{n\in\bZ}$ in $L^2(\bP)$ and the space $L^2(\rho)$ that maps $\xi (n)$
to $t^n$, $n\in\bZ$. Therefore, almost surely,
\begin{equation}\label{eq:annihil}
\sum_{n\in\bZ} \xi (n)\widehat{\phi}(n) = 0\,,
\end{equation}
that is, $F_\xi (\bar\phi)=0$. Since this holds for {\em any} $\phi$ as above, we conclude
that almost surely $\sigma (\xi)\subseteq \spt(\rho)$. \mbox{} \hfill $\Box$

\subsubsection{Proof of Part (B)}
One simply needs to read the same lines in the reverse direction. Suppose that
$I\subsetneqq\bT$ is an open arc such that almost surely $I\cap\sigma(\xi)=\emptyset$.
We take an arbitrary function $\phi\in C_0^\infty (I)$. Then, almost surely,
$F_\xi (\bar\phi)=0$, that is, \eqref{eq:annihil} holds. By the same linear isomorphism,
this yields~\eqref{eq:vanishing}, i.e., $\phi=0$ in $L^2(\rho)$. Since $\phi$ was
arbitrary, we conclude that $\spt(\rho)\cap I = \emptyset$. \mbox{} \hfill $\Box$

\subsection{Carleman spectrum}\label{subsect:Carleman}

In some instances, it is helpful to combine Theorem~\ref{thm:spectra} with another definition of the spectrum which is based on the analytic continuation and goes back at least
to Carleman.

Let $\xi\colon\bZ\to\bZ$ be a sequence of subexponential growth, i.e.,
$\displaystyle \limsup_{n\to\infty} |\xi (n)|^{1/n} = 1$. Put
\begin{align*}
F_\xi^+(z)&=\sum_{n\ge 0} \xi (n) z^n, \qquad
{\rm analytic\ in\ } \{|z|<1\}, \\
F_\xi^-(z)&=-\sum_{n\le -1} \xi (n) z^n, \qquad
{\rm analytic\ in\ } \{|z|> 1\}\cup \{\infty\}.
\end{align*}
Then the Carleman spectrum, $\sigma_C(\xi)$ is a minimal compact set $\sigma\subseteq\bT$
such that the function
\[
F_\xi (z) =
\begin{cases}
F_\xi^+(z), & |z|<1\, \\
F_\xi^- (z), & |z|>1
\end{cases}
\]
is analytic on $\widehat\bC\setminus\sigma$.

It is a classical fact of harmonic analysis that if the sequence $\xi$ has at most
polynomial growth (which is almost surely the case for realizations of a wide-sense
stationary process), then $\sigma_C(\xi)=\sigma(\xi)$, i.e., both definitions of the
spectrum of $\xi$ coincide (see, for example, \cite[VI.8]{Katz}, where the proof is given for bounded functions on $\bR$).

\subsection{Spectral sets of bounded sequences}
\label{subsect:Beurling}
It is worth to recall here that for bounded sequences $\xi$ there is another
equivalent definition of the spectrum that goes back to Beurling. The spectral
set (a.k.a. Beurling's spectrum or weak-star spectrum), $\sigma_B(\xi)$, of a
bounded sequence $\xi$ is the set of $t\in\bT$ such that the sequence $(t^n)_{n\in\bZ}$
belongs to the closure of the linear span of translates of $\xi$ in the weak-star
topology in $\ell^\infty(\bZ)$ as a dual to $\ell^1(\bZ)$.

The Hahn-Banach theorem easily gives an equivalent definition: $\sigma_B(\xi)$ is the set
of $t\in\bT$ such that whenever $\phi\in\ell^1(\bZ)$ and $\phi * \xi = 0$ identically,
necessarily $\sum_{n\in\bZ} \phi(n) \bar t^n = 0$.

It is not difficult to show that for any bounded sequence $\xi$,
$\sigma_B(\xi)=\sigma (\xi)$ (see for instance~\cite[VI.6]{Katz},
where the proof is given for bounded functions on $\bR$).

The advantage of the definition of the spectral set is that it can be
transferred to bounded functions on locally compact abelian groups.
This allows one without extra work to extend Theorem~\ref{thm:spectra}
to bounded stationary processes on arbitrary locally compact abelian
groups with a countable base of open sets.

\section{A random variation on the theme of Szeg\H{o} and Helson}\label{sect:Szego-Helson}

We call the wide-sense stationary random process $\xi$ {\em periodic},
if there exists a positive integer $N$ so that almost every realization
$(\xi (n))_{n\in\bZ}$ is $N$-periodic. In this case the spectrum of $\xi$
is contained in the set of roots of unity of order $N$.

\begin{theorem}\label{thm:wide-sense}
Let $X\subset\bC$ be a finite set, and let $\xi\colon\bZ\to X$ be a
wide-sense stationary process with the spectral measure $\rho$. Suppose that
$\spt(\rho)\ne\bT$. Then the process $\xi$ is periodic.
\end{theorem}

Theorem~\ref{thm:wide-sense} is an immediate corollary to the following theorem combined with
Theorem~\ref{thm:spectra}.

\subsection{Helson's theorem}\label{subsect:Helson}

In the following theorem which is essentially due to Helson \cite[Section~5.4]{Helson}
we have added the claim that the period $N$ does not depend on the
sequence $\xi$ but only on the spectrum and the set $X$. Helson's proof
used a compactness argument which doesn't give this and we have found a new
proof of his theorem which does.

\begin{theorem}[Helson]\label{thm:Helson}
Let $X\subset\bC$ be a finite set. Then any sequence $\xi\colon\bZ\to X$
with $\sigma(\xi)\ne\bT$ is $N$-periodic with $N$ depending only on
$X$ and on $\sigma (\xi)$. Moreover, $N$ as a function of
$\sigma(\xi)$ is non-decreasing.
\end{theorem}

%Note that the difference between Helson's version~\cite[Section~5.4]{Helson}
%and the one presented here is the claim that the period $N$ does not depend on
%the sequence $\xi$, but only its spectrum and the set $X$. Helson's proof used
%a compactness argument so we were unable to use it here.

\subsection{Two lemmas on polynomials on an arc of the circle}
\label{subsect:polynomials}

The proof of Theorem~\ref{thm:Helson} uses the following lemma:

\begin{lemma}\label{lemma:small-poly}
Given a closed arc $J\subsetneqq\bT$ and $\delta>0$, there exists a
polynomial $P$ such that $P(0)=1$ and $\| P \|_{C(J)}<\delta$.
\end{lemma}

\noindent{\em Proof:} By Runge's theorem, there exists a polynomial
$S$ so that $\bigl\| \tfrac1{z}-S \bigr\|_{C(J)}<\delta$. Put $P(z)=1-zS(z)$.
\mbox{} \hfill $\Box$

\medskip
Next, we will need a version of the Bernstein inequality for polynomials on
an arc which is due to V. S. Videnskii~\cite[Section~5.1.E19.c]{BE}.

\begin{lemma}\label{lemma-Bernstein}
Let $J\subsetneqq\bT$ be a closed arc. There exists a constant $K(J)$
such that for any polynomial $P$ of degree $n$, we have
$\| P'\|_{C(J)}\le K(J) n^2 \|P\|_{C(J)}$.
\end{lemma}

Note that any bound polynomial (or even subexponential)
in $n$ would suffice for our purposes.

\subsection{The $\delta$-prediction lemma}
\label{subsect: prediction}

The next lemma (which can be traced to Szeg\H{o}) is
the main ingredient in the proof of Theorem~\ref{thm:Helson}.
By $\|\xi\|_p$ we denote the weighted norm of $\xi$ defined as
\[
\| \xi \|_p^2 = \sum_{m\in\bZ} \Bigl( \frac{|\xi (m)|}{1+|m|^p} \Bigr)^2\,.
\]

\begin{lemma}\label{lemma:prediction}
Given $\delta, p, M >0$ and given an open arc $J$, $\bar J\subsetneqq\bT$,
there exist $n\in\bN$ and $q_0, \ldots , q_{n-1}\in\bC$ such that for any sequence
$\xi\colon \bZ\to \bC$ with $\|\xi\|_p\le M$ and $\sigma (\xi)\subset J$, we have
\[
\Bigl| \xi (n) + \sum_{k=0}^{n-1} q_k \xi (k) \Bigr|< \delta\,.
\]
\end{lemma}

\noindent{\em Proof}: It will be convenient to assume that
$J=\{e^{{\rm i}\theta}\colon |\theta|< \pi-\e \}$. Put
$J'=\{e^{{\rm i}\theta}\colon |\theta|< \pi-\tfrac12 \e \}$
and let $\phi\colon \bT\to [0, 1]$ be a $C^\infty$-function such
that $\phi=1$ on $J$ and $\phi=0$ on $\bT\setminus J'$.

By Lemma~\ref{lemma:small-poly}, there exists a polynomial $P$ such that
$P(0)=1$ and $\| P \|_{C(J')} \le \tfrac12$. Let $\ell\in\bN$ be a sufficiently
large number to be chosen later, $n=\ell\cdot\deg P$, and put
\[
Q(t) = t^{-n}P(t)^\ell = t^{-n} + \sum_{k=0}^{n-1} q_k t^{-k}\,.
\]
Then $\| Q \|_{C(J')} \le e^{-cn}$, whence by Lemma~\ref{lemma-Bernstein},
$\| Q^{(j)} \|_{C(J')} \le \bigl( K(J)n^2 \bigr)^j e^{-cn}$.

Next, we note that since the function $1-\phi$ vanishes on an open neighbourhood of the
spectrum of $\xi$, we have $F_\xi (Q) = F_\xi (Q\phi)$, and therefore,
\begin{align*}
\xi(n) + \sum_{k=0}^{n-1} q_k \xi(k) &= F_\xi (Q) \\
&= F_\xi (Q\phi) = \sum_{m\in\bZ} \xi (m) \widehat{Q\phi}(-m)\,.
\end{align*}
Hence, by the Cauchy-Schwarz inequality, we have
\begin{align*}
\Bigl| \xi(n) + \sum_{k=0}^{n-1} q_k \xi(k) \Bigr| &\le
\| \xi \|_p
\Bigl( \sum_{m\in\bZ} (1+|m|^p)^2 \bigl| \widehat{Q\phi}(-m) \bigr|^2 \Bigr)^{1/2} \\
&\lesssim_p M  \sum_{j=0}^p \| (Q\phi)^{(j)} \|_{L^2(\bT)} \\
&\lesssim_p  M  \sum_{j=0}^p \| (Q\phi)^{(j)} \|_{C(J')} \\
&\lesssim_{p, J} M  n^{2p} e^{-cn} \\
&\le \delta \qquad \text{for\ } n \ge n_0(p, J, M, \delta)\,,
\end{align*}
that completes the proof. \mbox{} \hfill $\Box$

\subsection{Proof of Theorem~\ref{thm:Helson}}
\label{subsect:proof-Helson}

We fix a finite set $X\subset\bC$ and put
$\delta_X = \inf\bigl\{|z-w|\colon z, w\in X, z\ne w \bigr\}$.
Applying the previous lemma with $\delta<\tfrac12 \delta_X$, we see that
the $n$-tuple $\langle \xi(0), \ldots , \xi(n-1) \rangle$ uniquely defines
the value of $\xi (n)$. Repeating this argument, we proceed with
$\xi (n+1), \xi (n+2), \ldots $. Since the sequences $(\xi (m))_{m\in\bZ}$
and $(\xi (-m))_{m\in\bZ}$ have the same spectrum, we proceed in the same way with
$\xi (-1), \xi (-2), \ldots\, $. Hence, the $n$-tuple
$\langle \xi(0), \ldots , \xi(n-1) \rangle$ uniquely determines the whole sequence
$\xi$.

Observe that among $|X|^n+1$ tuples
\[
\langle \xi(k), \xi (k+1), \ldots , \xi(k + n-1) \rangle, \quad
0\le k \le |X|^n,
\]
at least two must coincide. If $0\le k_1 < k_2 \le |X|^n$ are the corresponding indices,
after a minute reflection we conclude that the sequence $\xi$ is periodic with
period $N=k_2-k_1 \le |X|^n$. \mbox{} \hfill $\Box$

\section{More random variations on Szeg\H{o}'s theme}
\label{sect:2nd-variation}
Here we will give another version of our main result, Theorem~\ref{thm:wide-sense}.
Now we will assume that the random process $\xi\colon \bZ\to X$ is stationary in
the usual sense and replace our previous assumption that the spectrum of $\xi$ is
not all of $\bT$ by a much weaker condition imposed on the spectral measure $\rho$ of $\xi$.
We will also allow $X$ to be a uniformly discrete subset of $\bC$.

\subsection{Condition $(\Theta)$}
\label{subsect:condition_Theta}

Given a positive measure $\rho$ on $\bT$, we put
\[
e_n(\rho) = \operatorname{dist}_{L^2(\rho)} \bigl(\done, \mathcal P_n^0 \bigr),
\]
where $\mathcal P_n^0 \subset\bC[z]$ is the set of algebraic polynomials of degree at most
$n$ vanishing at the origin. We say that {\em the measure $\rho$ satisfies condition $(\Theta)$} if
\[
\sum_{n\ge 1} e_n(\rho)^2 < \infty\,.
\]
Recall that by the classical Szeg\H{o} theorem,
\[
\lim_{n\to\infty} e_n(\rho) = 0 \ \Longleftrightarrow \int_{\bT} \log \rho'\, {\rm d}m = -\infty\,,
\]
where $m$ is the Lebesgue measure on $\bT$ and
$\displaystyle \rho' = \frac{{\rm d}\rho}{{\rm d}m}$.
Note that by Lemma~\ref{lemma:small-poly}, $e_n(\rho)=O(e^{-cn})$ provided
that $\spt(\rho)\ne\bT$, that is, condition $(\Theta)$ is essentially weaker than
the condition that the support of $\rho$ is not all of $\bT$. In Section~\ref{sect:Theta} we
will show that condition $(\Theta)$ forbids $\rho$ to have certain exponentially deep zeroes.

In spite of the omnipresence of Szeg\H{o}'s theorem~\cite{GS} (see
also~\cite{Bingham, DIK, Simon} for recent development) we are not aware of
any results that in the case when the logarithmic integral of $\rho'$ is divergent would relate the rate of decay of the sequence $(e_n(\rho))$ with the properties of the measure $\rho$ (the only exceptions are several results that relate the gaps in $\spt(\rho)$
with the exponential rate of decay of $e_n(\rho)$'s).

\subsection{Stationary processes on $\bZ$ with uniformly discrete set of values}
\label{subsect:unif_discrete}

Recall that the set $X\subset\bC$ is called {\em uniformly discrete} if
\[
\delta_X \stackrel{\rm def}=\inf\bigl\{|z-w|\colon z, w\in X, z\ne w \bigr\}>0\,.
\]

\begin{theorem}\label{thm:stationary}\mbox{}
Suppose $X\subset\bC$ is a uniformly discrete set of values
and $\xi\colon\bZ\to X$ is a stationary process with the spectral measure $\rho$ satisfying
condition $(\Theta)$. Then almost every realization $(\xi (n))_{n\in\bZ}$ of the process $\xi$ is periodic.
\end{theorem}

\begin{corollary}\label{cor:unif_discrete}
Suppose that the stationary ergodic process $\xi$ on $\bZ$ attains values
in a uniformly discrete subset of $\bC$ and that its spectral measure $\rho$
satisfies condition $(\Theta)$. Then the process $\xi$ is periodic.
\end{corollary}

This was proven in~\cite{BNS} under a stronger assumption that $\spt(\rho)\ne\bT$. The proof given in~\cite{BNS} used a somewhat different approach.

\subsection{Proof of Theorem~\ref{thm:stationary}}
\label{subsect:Proof_part(A)}

Applying the linear isomorphism between the closure of the linear span of
$(\xi (n))_{n\in\bZ}$ in $L^2(\bP)$ and the space $L^2(\rho)$ that maps $\xi (n)$ to $t^n$, we see that, given $N\in\bN$, there exist $q_0, \ldots , q_{N-1}\in\bC$ so that
\[
\bE\Bigl[ \bigl| \xi (N) + \sum_{k=0}^{N-1} q_k \xi (k) \bigr|^2 \Bigr] = e_N(\rho)^2\,,
\]
whence
\[
\bP\Bigl\{ \bigl| \xi (N) + \sum_{k=0}^{N-1} q_k \xi (k) \bigr| \ge \frac12 \delta_X \Bigr\}
\le \frac4{\delta_X^2}\, e_N(\rho)^2\,.
\]
Thus, with probability at least $1-\tfrac4{\delta_X^2}e_N(\rho)^2$, the term $\xi (N)$ is determined
by the values of the terms $\xi (0), \ldots , \xi (N-1)$. Proceeding to $\xi (N+1), \xi(N+2), \ldots$
and then to $\xi (-1), \xi (-2), \ldots $, we see that the whole sequence $(\xi (n))$
is determined by the values $\xi (0), \ldots , \xi (N-1)$ with probability at least
\[
1 - \frac8{\delta_X^2}\, \sum_{n\ge N} e_n(\rho)^2\,.
\]
Recalling that the series of $e_n(\rho)^2$ converges (and that the set $X$ is countable), we
conclude that, given $\e>0$, there is a countable part $\Omega'$ of our probability space
$\Omega$ with $\bP\{\Omega'\}>1-\e$. Hence, we may assume that the whole probability space $\Omega$ is countable.

Since the process $\xi$ is a stationary one, there exists a measure preserving
transformation $\tau$ of $(\Omega, \bP)$. The only possibility for this is that
$\Omega=\bigcup_j \Omega_j$, each $\Omega_j$ consists of finitely many atoms of
equal probability and is $\tau$-invariant. That is, $\tau$ acts as a permutation
of elements of $\Omega_j$. Therefore, for each $j$, the restriction $\tau|_{\Omega_j}$ is periodic.
\mbox{}\phantom{A} \hfill $\Box$

\subsection{Integer-valued stationary processes}
Here is another variation on the same theme taken from~\cite{BNS}:
\begin{theorem}\label{thm:integer-valued}
Suppose $\xi\colon \bZ\to\bZ$ is a stationary process with $\spt(\rho)\ne\bT$.
Then the process $\xi$ is periodic.
\end{theorem}

\subsection{Proof of Theorem~\ref{thm:integer-valued}}
\label{subsect:Proof_part(B)}

By Theorem~\ref{thm:stationary}, almost every realization $\xi (n)$ is periodic with some period $N$.
We need to show that the period $N$ is non-random. For this, it suffices to see that
it is bounded.

We have
\begin{align*}
F_\xi^+(z) &\stackrel{\rm def}= \sum_{n\ge 0} \xi (n) z^n \\
&= (1-z^N)^{-1}\, \sum_{n=0}^{N-1} \xi (n) z^n \\
&=\frac{T(z)}{1-z^N}\,,
\end{align*}
where $T$ is a polynomial with integer coefficients of degree at most $N-1$.
Since the ring $\bZ[z]$ of polynomials over the integers is a unique factorization domain, we can
write
\[
F_\xi^+(z) = \frac{P(z)}{Q(z)}\,,
\]
where $P, Q \in \bZ[z]$, have no common factors, and $Q$ is monic and divides
$z^N-1$.

First, we claim that $P$ and $Q$ have no common zeroes. Indeed,
we decompose $P=P_1 \cdot \ldots \cdot P_\ell$ and $Q=Q_1\cdot \ldots \cdot Q_m$
in the product of monic polynomials irreducible on $\bZ$. By the Gauss lemma, these
polynomials are also irreducible on $\bQ$. Since $\bQ$ is a field, $\bQ[z]$ is
a principal ideal domain, that is, if $P_i\ne Q_j$ are irreducible on $\bQ$, then there
are $R, S\in \bQ[z]$ so that $P_i R + Q_j S = 1$. Thus, $P_i$ and $Q_j$ have no
common zeroes. Hence, $P$ and $Q$ have no common zeroes.

Next, let
\[
\Phi_n(z) = \prod_{k\colon {\rm gcd}(k, n)=1} ( z - e(k/n))
\]
be the cyclotomic polynomial of degree $n$. It belongs to $\bZ[z]$ and is irreducible~\cite[\S~8.4]{vdW}. Therefore, we can factor
\[
Q = \prod_{1\le k \le r} \Phi_{n(k)}\,.
\]
Since the rational function $F_\xi^+$ does not have singularities on a fixed arc of $\bT$,
there is a (non-random) integer $n^*$ so that $n(k)\le n^*$. Thus, there is a (non-random) integer $N^*$ so that
\[
{\rm zeroes}(Q) \subset \{t\colon t^{N^*} = 1 \}\,.
\]
That is, the period of the random sequence $(\xi (n))_{n\in\bZ}$ is bounded by $N^*$.
\mbox{} \hfill $\Box$

\section[Measures satisfying condition \footnotesize{$(\Theta)$}]{Measures satisfying condition $(\Theta)$}
\label{sect:Theta}

The following theorem yields that if a positive measure $\rho$ on $\bT$
has a deep exponential root, then it must satisfy condition $(\Theta)$,
and if the exponential root is not deep enough, then $\rho$ may not satisfy
condition $(\Theta)$.

As before, we denote by $m$ the normalized
Lebesgue measure on $\bT$.

\begin{theorem}\label{thm:exp-root}
Let $\rho$ be a positive measure on $\bT$ and let
$\beta$ be a positive parameter.

\smallskip\noindent{\rm (A)}
There exists a positive numerical constant $K$ such that if
\[
\int_{-\pi}^\pi \exp\Bigl( \frac{\beta}{|\theta|} \Bigr)\, {\rm d}\rho(e^{{\rm i}\theta})
< \infty\,,
\]
then $e_n(\rho)\lesssim n^{-K \beta}$.

\smallskip\noindent{\rm (B)}
Suppose that $\displaystyle {\rm d}\rho \gtrsim
\exp\bigl( -\frac{\beta}{|\theta|} \bigr)\, {\rm d}m$ everywhere on $\bT$.
Then $e_n(\rho) \gtrsim_\beta n^{-\beta/(2\pi)}$.
\end{theorem}

We will deduce Theorem~\ref{thm:exp-root} from a more general estimate
which in spite of its crudeness might be of an independent interest.
We let $W\colon\bT\to [1, +\infty]$ be an arbitrary measurable function, and
for $A>0$ we put $W_A=\min (W, e^A)$. By
\[
P_r (t) = \frac{1-r^2}{|1-rt|^2}, \quad 0\le r<1, \ t\in\bT,
\]
we denote the Poisson kernel for the unit disk.

\begin{theorem}\label{thm:Theta} Let $\rho$ be a positive measure on
$\bT$, let $\beta, M$ be positive parameters, and let $A\ge A_0(\beta, M) >1 $.

\smallskip\noindent{\rm (A)}
Suppose that
\begin{equation}\label{eq:1}
\int_{\bT} W^\beta\, {\rm d}\rho \le C <\infty\,,
\end{equation}
and that
\begin{equation}\label{eq:2}
\bigl( \log W_A \bigr) * P_{1-A^{-1}} \le M \log W
\end{equation}
everywhere on $\bT$. Then, for $\displaystyle n\ge \frac{A^2 \beta}{2M}$,
\[
e_n(\rho) \le \sqrt{2C+\tfrac12 \rho(\bT)}\, \exp\Bigl( -\frac{\beta}{2M}\, \int_{\bT} \log W_A\, {\rm d} m\Bigr)\,.
\]

\smallskip\noindent{\rm (B)} Suppose that ${\rm d}\rho \gtrsim W^{-\beta}\, {\rm d}m$ everywhere on $\bT$. Then, for $n \cdot m\{ W>e^A\} \lesssim 1$,
\[
e_n(\rho) \gtrsim \exp\Bigl( -\frac{\beta}2\, \int_{\bT} \log W_A\, {\rm d} m\Bigr)\,.
\]
\end{theorem}

\subsection{Proof of Theorem~\ref{thm:Theta}, Part (A)}

Let $F_A$ be an outer function in the unit disk with the boundary values
$|F_A|^2 = W_A^{\beta/M}$ on $\bT$ (this means that
\[
\log |F_A(rt)| = \bigl( \frac{\beta}{2M} \log W_A * P_r \bigr)(t)
\]
everywhere in the unit disk $\bD$). Then, by \eqref{eq:2},
\[
\int_\bT \bigl| F_A\bigl( (1-A^{-1})t \bigr) \bigr|^2\, {\rm d}\rho(t)
\le \int_\bT W^\beta\, {\rm d}\rho \le C\,.
\]
Furthermore,
$ |F_A| \le e^{A\beta/(2M)} $
everywhere in $\bD$, and
\begin{align*}
|F_A(0)| &=
\exp\Bigl( \int_\bT \log|F_A|\, {\rm d}m \Bigr) \\
&= \exp\Bigl( \frac{\beta}{2M} \int_\bT \log W_A \,{\rm d}m \Bigr)\,.
\end{align*}

Let
\[
F_A\bigl( (1-A^{-1})z \bigr) = \sum_{k\ge 0} f_k z^k\,.
\]
Then, by Cauchy's inequalities,
\[
|f_k| \le \bigl( 1-A^{-1} \bigr)^k\, e^{A\beta/(2M)}\,,
\]
and for $n\ge A^2\beta/M$ we get
\begin{align*}
\sum_{k>n} |f_k| &\le \bigl( 1-A^{-1} \bigr)^n\, e^{A\beta/(2M)} \cdot A \\
&\le A \exp\Bigl( \frac{A\beta}{2M} - \frac{n}A \Bigr) \\
&\le A \exp\Bigl( - \frac{A\beta}{2M} \Bigr) \\
&< \frac12\,, \quad \text{provided\ that\ } A\ge A_0(\beta, M).
\end{align*}

Put
\[
P_A(z) = \sum_{k=0}^n f_k z^k\,.
\]
Then
\[
\int_{\bT} |P_A|^2\, {\rm d}\rho \le 2\int_{\bT} |F_A|^2{\rm d}\rho + \frac12 \rho(\bT)
\le 2C+\frac12 \rho(\bT)\,,
\]
while
\[
|P_A(0)| = |F_A(0)|
\ge \exp\Bigl( \frac{\beta}{2M}\, \int_\bT \log W_A\, {\rm d}m \Bigr)\,.
\]
Therefore,
\[
e_n(\rho) \le \sqrt{2C+\tfrac12 \rho(\bT)}\, \exp\Bigl( -\frac{\beta}{2M}\, \int_{\bT} \log W_A\, {\rm d} m\Bigr)\,,
\]
proving part (A). \mbox{} \hfill $\Box$

\subsection{Proof of Theorem~\ref{thm:Theta}, Part (B)}
To simplify notation, we denote here by $C$, $c$  various positive numerical constant whose values are not
important for our purposes. The values of these constants may vary from line to line.

Let $P$  be a polynomial of degree $n$  such that $P(0)=1$  and
\[
\int_\bT |P|^2\, {\rm d}\rho = e_n(\rho)^2\,.
\]
Then
\begin{align*}
0 = \log |P(0)| &\le \int_\bT \log |P|\, {\rm d}m  \\
&= \frac12\, \int_\bT \log\bigl( |P|^2 W_A^{-\beta} \bigr)\, {\rm d}m
+\frac{\beta}2\, \int_\bT \log W_A\, {\rm d}m\,.
\end{align*}
Furthermore,
\begin{align*}
\int_\bT |P|^2 W_A^{-\beta}\, {\rm d}m &=
\int_{\{W\le e^A\}}  |P|^2 W_A^{-\beta}\, {\rm d}m
+ e^{-\beta A}\, \int_{\{W>e^A\}} |P|^2\, {\rm d}m \\
&\le \int_{\{W\le e^A\}}  |P|^2 W_A^{-\beta}\, {\rm d}m +
e^{-\beta A}\, \int_\bT  |P|^2\, {\rm d}m\,,
\end{align*}
and by Nazarov's version of Tur\'an's lemma~\cite{Nazarov}
(for the reader's convenience, we will recall it after %we will complete
the proof):
\begin{align*}
\int_\bT |P|^2\, {\rm d}m &\le e^{Cn\cdot m\{W>e^A\}}\, \int_{{W\le e^A}} |P|^2\, {\rm d}m \\
&\le C\, \int_{{W\le e^A}} |P|^2\, {\rm d}m\,,
\end{align*}
provided that $ n\cdot m\{W>e^A\} \le C$. Thus,
\begin{align*}
\int_\bT |P|^2 W_A^{-\beta}\, {\rm d}m &\le
\int_{\{W\le e^A\}} |P|^2 W_A^{-\beta}\, {\rm d}m
+ Ce^{-\beta A}\, \int_{\{W\le e^A\}} |P|^2\, {\rm d}m \\
%&\le \int_{\{W\le e^A\}} |P|^2 W_A^{-\beta}\, {\rm d}m
%+ C\, \int_{\{W\le e^A\}} |P|^2 W_A^{-\beta}\, {\rm d}m \\
&\lesssim  \int_\bT |P|^2 W_A^{-\beta}\, {\rm d}m \\
%&\lesssim \int_\bT |P|^2 W^{-\beta}\, {\rm d}m \\
&\lesssim \int_\bT |P|^2\, {\rm d}\rho \\
&= e_n(\rho)^2 \quad \text{by\ the\ choice\ of\ the\ polynomial\ } P.
\end{align*}

Therefore,
\[
0 \le \frac12\log\bigl( e_n(\rho)^2 \bigr) + C + \frac{\beta}2\, \int_\bT \log W_A\, {\rm d}m\,,
\]
whence,
\[
e_n(\rho) \gtrsim \exp\Bigl( -\frac{\beta}2\, \int_\bT \log W_A\, {\rm d}m \Bigr)
\]
completing the proof. \mbox{} \hfill $\Box$

\medskip Here is the result of Nazarov used in the proof:
\begin{theorem}[Nazarov]\label{thm:Nazarov}
There exists a positive numerical constant $A$ such that for any set $\Lambda\subset\bZ$ of cardinality
$n+1$, any trigonometric polynomial
\[
p(t) = \sum_{\la\in\Lambda} c_\lambda t^\lambda\,, \quad (c_\la)_{\la\in\Lambda}\subset\bC, \ t\in\bT\,,
\]
and any measurable set $E\subset\bT$ of Lebesgue measure $m(E)\ge \tfrac13$,
\[
\int_\bT |p|^2\, {\rm d}m \le e^{An\cdot m(\bT\setminus E)}\, \int_E |p|^2\, {\rm d}m\,.
\]
\end{theorem}

We have used this theorem with $\Lambda=\{0, 1, \ldots , n\}$.

\subsection{Proof of Theorem~\ref{thm:exp-root}}
We apply Theorem~\ref{thm:Theta} with
\[
W(e^{{\rm i}\theta}) = \exp\bigl( \frac1{|\theta|} \bigr), \qquad -\pi \le \theta \le \pi\,.
\]
Then
\[
\Bigl| \int_{\bT} \log W_A\, {\rm d}m - \frac1{\pi}\, \log A \Bigr| \le C\,,
\]
and
$ m\{W>e^A\} = (\pi A)^{-1}$. Thus, choosing $A\simeq n$, we get
\[
e_n(\rho) \gtrsim \exp\Bigl(- \frac{\beta}{2\pi} \log A \Bigr)
\gtrsim_\beta n ^{-\beta/(2\pi)}\,,
\]
proving the lower bound.

To get the upper bound we need to check condition~\eqref{eq:2}. This is straightforward but somewhat long.
To simplify our writing, we put
\[
w_A(\theta)=\min\bigl( \frac1{|\theta|}, A \bigr)
\quad \text{and} \quad p_r(\theta) = P_r(e^{{\rm i}\theta})\,,
\quad -\pi \le \theta \le \pi\,.
\]
In this notation, we will show that
\begin{equation}\label{eq:2'}
\int_{-\pi}^\pi w_A(\theta-\phi) p_{1-A^{-1}}(\phi)\, {\rm d}\phi
\lesssim w_A(\theta)\,.
\end{equation}
Since the function $w_A$ is even, we may assume that $0\le\theta\le\pi$.

Using the estimate of the Poisson kernel
\[
p_r(\phi) \lesssim (1-r)\min\bigl(\phi^{-2}, (1-r)^{-2} \bigr),
\]
we get
\begin{align*}
\int_{-\pi}^\pi w_A(\theta-\phi) p_{1-A^{-1}}(\phi)\, {\rm d}\phi
&\lesssim
\frac1{A}\, \int_{-\pi}^\pi w_A(\theta-\phi) \min\bigl( \phi^{-2}, A^2 \bigr)\, {\rm d}\phi
\\
&=\frac1{A}\, \int_0^\pi
\bigl( w_A(\theta-\phi) + w_A(\theta+\phi) \bigr)
\min\bigl( \phi^{-2}, A^2 \bigr)\, {\rm d}\phi \\
&\le \frac2{A}\, \int_0^\pi
w_A(\theta-\phi) \min\bigl( \phi^{-2}, A^2 \bigr)\, {\rm d}\phi \\
&\lesssim A \int_0^{1/A} w_A(\theta-\phi)\, {\rm d}\phi
+ \frac1{A}\, \sum_{j\ge 0} \int_{2^j/A}^{2^{j+1}/A} w_A(\theta-\phi) \phi^{-2}\, {\rm d}\phi\,,
\end{align*}
It remains to estimate the integrals. We have
\begin{align*}
A \int_0^{1/A} w_A(\theta-\phi)\, {\rm d}\phi
&= A \int_0^{1/A} \bigl( w_A(\theta) + O(1) \min(\theta^{-2}, A^2)\phi\bigr)\, {\rm d}\phi \\
&= w_A(\theta) + O(1)A^{-1}\min(\theta^{-2}, A^2) \\
&\lesssim w_A(\theta)\,.
\end{align*}
Next,
\[
\frac1{A}\, \sum_{j\ge 0} \int_{2^j/A}^{2^{j+1}/A} w_A(\theta-\phi) \phi^{-2}\, {\rm d}\phi
\le A \sum_{j\ge 0} 2^{-2j} \int_{2^j/A}^{2^{j+1}/A} w_A(\theta-\phi)\, {\rm d}\phi\,.
\]
Let $k$ be the minimal non-negative integer such that $\theta\le 2^k/A$. Then
\begin{align*}
A \sum_{j\ge k+1} 2^{-2j} \int_{2^j/A}^{2^{j+1}/A} w_A(\theta-\phi)\, {\rm d}\phi
&=
A \sum_{j\ge k+1} 2^{-2j} \int_{2^j/A}^{2^{j+1}/A} \frac{{\rm d}\phi}{\phi-\theta} \\
&\lesssim
A \sum_{j\ge k+1} 2^{-2j} \int_{2^j/A}^{2^{j+1}/A} \frac{{\rm d}\phi}{\phi}\\
&\lesssim
A \sum_{j\ge k+1} 2^{-2j} \\
%&\lesssim A\, 2^{-2k} \\
&\lesssim 2^{-k} w_A(\theta)\,,
\end{align*}
and, for $k\ge 2$,
\begin{align*}
A \sum_{0\le j\le k-2} 2^{-2j} \int_{2^j/A}^{2^{j+1}/A} w_A(\theta-\phi)\, {\rm d}\phi
&= A \sum_{0\le j\le k-2} 2^{-2j} \int_{2^j/A}^{2^{j+1}/A} \frac{{\rm d}\phi}{\theta-\phi} \\
&\lesssim \frac1{\theta}\, A\, \sum_{0\le j \le k-2} 2^{-2j}\, \frac{2^j}A \\
&\lesssim \frac1{\theta}\,.
\end{align*}
At last, we are left with two terms corresponding to $j=k$ and $j=k-1$ (the second term disappears when $k=0$).
Each of them is bounded by
\[
\frac{A}{2^{2k}}\cdot A \cdot\frac{2^k}A = \frac{A}{2^k}
\lesssim w_A(\theta).
\]
This completes the proof of estimate~\eqref{eq:2'}, and hence of Theorem~\ref{thm:exp-root}.
\mbox{} \hfill $\Box$

\section{Unimodular processes on $\bZ$}
\label{sect:unimodular}

Here we are interested in spectral properties of unimodular stationary processes
$\xi\colon\bZ\to\bT$.

\subsection{An example}\label{subsect:unimod_example}
We start with a simple example that shows that {\em any} probability measure
on $\bT$ can be realized as the spectral measure of a unimodular
stationary process on $\bZ$.

Consider the product of two circles $\bT\times \bT$, the first equipped with a
probability measure $\rho$ and the second with the normalized Lebesgue measure,
and let $P$ be the product of these measures on $\bT\times\bT$. Let $\tau$ be the
measure preserving transformation which is the identity on the first circle, and
rotates the second by the first, i.e.,
$ \tau (s, t) = (s, st) $, $(s, t)\in\bT\times\bT$.
Consider the process $\xi$ defined by the function $f(s, t)=t$ and the transformation
$\tau$, i.e., put $\xi (n) = f(\tau^n) = t\cdot s^n$, $n\in\bZ$. It remains to note that
the covariance function of this process is the Fourier transform of the measure $\rho$:
\[
r(m) = \bE[\xi(0)\bar\xi(m)] = \int_{\bT} {\bar s}\,^m\, {\rm d}\rho (s)
= \widehat{\rho}(m)\,, \qquad m\in\bZ\,.
\]

\subsection{A theorem}\label{subsect:unimod_theorem}

Curiously enough, it appears that if the spectrum of the unimodular stationary process
is contained in an arc of length less than $\pi$, then the only possibility is given
by the example we gave.

\begin{theorem}\label{thm:unimodular}
Let $\xi\colon\bZ\to\bT$ be a unimodular wide-sense stationary process
with the spectral measure $\rho$. Suppose that $\spt(\rho)$ is contained in an
arc of length less than $\pi$. Then almost every realization of $\xi$ has the
form
\begin{equation}\label{eq:unimod}
\xi(n) = t \cdot  s^n, \qquad n\in\bZ\,,
\end{equation}
with some (random) $t,s\in\bT$.
\end{theorem}

\begin{corollary}\label{cor:unimod}
The only unimodular stationary ergodic process $\xi$ on $\bZ$ with the support
of the spectral measure contained in an arc of length less than $\pi$ is
the process defined by~\eqref{eq:unimod} for some constant $s\in\bT$
and for $t\in\bT$ which is
a random variable with either a uniform distribution on the circle if $s$ is an irrational rotation, or a uniform distribution on a cyclic group
if $s$ is a rational rotation.
\end{corollary}

Indeed, the ergodicity of $\xi$ yields that $s\in\bT$ is a constant. Furthermore, the distribution of $t\in\bT$
is invariant under multiplication by $s$ and ergodic. \mbox{} \hfill $\Box$

\medskip
It would be very interesting to reveal other restrictions on the spectral measures
of unimodular ergodic processes on $\bZ$.

\subsubsection{Eremenko--Ostrovskii theorem}

The proof of Theorem~\ref{thm:unimodular}
is based on the following deterministic result of Eremenko and Ostrovskii~\cite[Theorem~$1'$]{EO}:

\begin{theorem}[Eremenko-Ostrovskii]\label{thm:EO}
Suppose that the Taylor series
\[
F(z) = \sum_{n\ge 0} f_n z^n
\]
with unimodular coefficients $(f_n)_{n\ge 0}\subset\bT$ has an analytic continuation on $\bar\bC\setminus J$, where
$J\subset\bT$ is an arc of length less than $\pi$. Then $f_n=t\cdot s^n$, $n\ge 0$, with $t, s\in\bT$.
\end{theorem}

\subsubsection{Proof of Theorem~\ref{thm:unimodular}}
\label{subsect:unimod-proof}

Let $(\xi (n))_{n\in\bZ}$ be a realization of the process $\xi$.
By Theorem~\ref{thm:spectra}, the Fourier-Carleman transform $F_\xi$ is analytic on
$\bar\bC\setminus J$ where $J\subset\bT$ is an arc of length less than $\pi$.
Then, applying Theorem~\ref{thm:EO} first
to $F_\xi^+$ and then to $F_\xi^-$, we conclude that~\eqref{eq:unimod} holds
with some $t,s\in\bT$. \mbox{} \hfill $\Box$

\bigskip
\medskip

{\noindent A.B.:
Aix-Marseille Universit\'e, CNRS, Centrale Marseille, I2M, 13453  Marseille,
France,
\newline {\tt borichev@cmi.univ-mrs.fr}

\smallskip
\noindent M.S.:
 School of Mathematics, Tel Aviv University, Tel Aviv, 69978 Israel,
\newline  {\tt sodin@post.tau.ac.il}

\smallskip
\noindent B.W.:
Institute of Mathematics, Hebrew University of Jerusalem,
Jerusalem, 91904 Israel,
\newline {\tt weiss@math.huji.ac.il}
}

\end{document}